\documentclass[12pt]{article}

\usepackage{amsmath,amssymb,amsfonts,theorem,makeidx,latexsym,epsfig,subfigure}
\usepackage{color}

\newtheorem{defn}{Definition}[section]

\newtheorem{lemma}[defn]{Lemma}

{\theorembodyfont{\rmfamily}

\newtheorem{ex}[defn]{Example}}

\newtheorem{thm}[defn]{Theorem}

\newtheorem{prop}[defn]{Proposition}

\newtheorem{cor}[defn]{Corollary}

\numberwithin{equation}{section}

\newcommand{\h}{{\cal H}}

\newcommand{\lk}{\lambda_k}

\newcommand{\ltr}{ L^2(\mathbb R) }

\newcommand{\mn}{\mathbb N}

\newcommand{\mr}{\mathbb R}

\newcommand{\mz}{\mathbb Z}

\newcommand{\mts}{ \{E_{mb}T_{na}g \}_{m,n \in \mz}}

\def\bp{{\noindent\bf Proof. \ }}

\def\ep{\hfill$\square$\par\bigskip}

\def\bqs{\begin{equation}}

\def\eqs{\tag*{$\square$}\end{equation}\par\bigskip}

\def\la{\langle}

\def\ra{\rangle}

\def\ga{\gamma}

\def\ftk{\{f_k\}_{k=1}^\infty}

\def\etk{\{e_k\}_{k=1}^\infty}

\def\suk{\sum_{k=1}^\infty}

\def\supp{\text{supp}}

\def\vn{\vspace{.1in}\noindent}

\def\bop{\begin{op}\rm}

\def\eop{\end{op}}

\def\bee{\begin{eqnarray}}

\def\ene{\end{eqnarray}}

\def\bes{\begin{eqnarray*}}

\def\ens{\end{eqnarray*}}

\def\bei{\begin{itemize}}

\def\eni{\end{itemize}}

\def\bt{\begin{thm}}

\def\et{\end{thm}}

\def\bc{\begin{cor}}

\def\ec{\end{cor}}

\def\bpr{\begin{prop}}

\def\epr{\end{prop}}

\def\bl{\begin{lemma}}

\def\el{\end{lemma}}

\def\bd{\begin{defn}}

\def\ed{\end{defn}}

\def\bex{\begin{ex}}

\def\enx{\end{ex}}

\def\bfi{\begin{fig}}

\def\efi{\end{fig}}

\def\lk{\lambda_k}

\def\mno{\mn_0}

\def\lkn{\{\lambda_k\}_{k=1}^\infty}

\title{The mystery of Carleson frames}

\date{\today}

\author{Ole Christensen, Marzieh Hasannasab, Friedrich Philipp, Diana Stoeva }

\begin{document}

\maketitle

\begin{abstract} In 2016 Aldroubi et al. constructed the
	first class of frames having the form $\{T^k \varphi\}_{k=0}^\infty$ for a bounded linear operator
	on the underlying Hilbert space. In this paper
	we show that a subclass of these frames has a number of additional remarkable features that have not  been
	identified for any
	other frames in the literature. Most importantly, 
	the
	subfamily obtained by selecting each $N$th element from the frame is itself a frame, regardless of the choice of $N\in \mn.$ Furthermore,  the frame property is kept upon removal of an arbitrarily finite number of elements.

\end{abstract}

\section{Introduction}

Let $\mathbb{D}$ denote the open unit disc in the complex plane. Recall that a sequence $\{\lambda_k\}_{k=1}^\infty\subset \mathbb{D}$  is said
to satisfy  the Carleson condition if

\bee\label{carleson}\displaystyle\inf_{n\in \mn} \prod_{k\neq n}\frac{|\lk-\lambda_n|}{|1-\overline{\lk}\,\lambda_n|}>0.\ene

 In \cite{A1} Aldroubi et al.  proved the following result:
\bt \label{19408bf} Assume that $\{\lambda_k\}_{k=1}^\infty\subset \mathbb{D}$  satisfies the Carleson condition, 
and let $\etk$ denote an orthonormal basis for a separable Hilbert space $\h.$  Consider the
bounded linear operator $T: \h \to \h$
given by $Te_k:= \lk e_k, \, k\in \mn,$ and let $\{m_k\}_{k=1}^\infty$
be any scalar sequence such that $0< C_1 \le | m_k| \le C_2 < \infty.$  Then, with
\bee \label{232111a} \varphi:= \suk m_k \sqrt{1-|\lk |^2} \, e_k,\ene the family
$\{T^{k}\varphi\}_{k=0~}^\infty $ is a frame for
$\h.$ \et

For the sake of easy reference we will call the frames
in Theorem \ref{19408bf} for {\it Carleson frames.} 
The construction in Theorem \ref{19408bf} provides
the first examples of  frames that are given as the orbit
of a bounded operator.  Except for the Riesz bases (see
\cite{olemmaarzieh})
this class still forms the only explicitly known example
of such frames.  The purpose of the current paper is to
show that under a  weak extra condition, such frames have a number of 
important features that have not been identified for any other frames in the literature. Most importantly, the subsequence of 
$\{T^{k}\varphi\}_{k=0~}^\infty $ formed by selecting
each $N$th element of the sequence is itself a frame, regardless of the choice of $N\in \mn.$

Our frame results are based on a number of manipulations 
on the Carleson condition. For this reason we first
collect what is known about the Carleson condition in
the literature.

\begin{prop}\label{1406a}
Given a sequence $\lkn\subset\mathbb{D}$ of distinct numbers, the following
hold:

\bei
\item[(i)] {\bf \cite{duren70}} If
	\bee \label{1606c} \exists\, c\in(0,1) \mbox{ such that } \frac{1-|\lambda_{k+1}|}{1-|\lambda_{k}|}\leq c<1, \quad \forall k\in\mn,\ene
	then $\lkn$ satisfies the Carleson condition.
\item[(ii)]	{\bf \cite{duren70}} If $\lkn$ is positive and increasing,  then the condition \eqref{1606c} is also necessary for $\lkn$ to satisfy the Carleson condition.
\item[(iii)] {\bf \cite{shapiro}} If $\lkn$ satisfies the Carleson condition, then $\lim_{k\to \infty}|\lk| =1.$
\item[(iv)] {\bf \cite{CMPP}} If  there is some $n\in\mn$ such that $\{\lk\}_{k\geq n}$ satisfies the Carleson condition, then also $\lkn$ satisfies the Carleson condition.
    \eni
\end{prop}


\section{Frame subsequences of the Carleson frames}

In this section we will consider a subclass of the Carleson frames in Theorem \ref{19408bf}; more precisely, we will
add a weak assumption on the sequence $\{\lambda_k\}_{k=1}^\infty\subset \mathbb{D}.$  The resulting class
of Carleson frames have the  remarkable properties mentioned
in the abstract.

\bt \label{19408b} Let $\{\lambda_k\}_{k=1}^\infty\subset \mathbb{D}$ and
assume that $\{|\lambda_k|\}_{k=1}^\infty$ is a strictly
increasing sequence  which satisfies the Carleson condition. Let $\etk$ denote an orthonormal basis for a separable Hilbert space $\h, $ consider the
bounded linear operator $T: \h \to \h$
given by $Te_k:= \lk e_k, \, k\in \mn,$ and
 let $\{m_k\}_{k=1}^\infty$
be any scalar sequence such that $0< C_1 \le | m_k| \le C_2 < \infty.$ Then, with
\bee \label{230412p} \varphi:= \suk m_k \sqrt{1-|\lk |^2} \, e_k,\ene the
following hold:

\bei\item[(i)] The  family
$\{T^{Nk}\varphi\}_{k=0}^\infty $ is a frame for $\h$ for
any $N\in \mn.$
 \item[(ii)] If $\lambda_1 \neq 0,$ the family $\{T^{Nk} \varphi\}_{k=K}^\infty$ is a frame for $\h$ for any $N,K \in \mn.$ 
\item[(iii)] If $\lambda_1\neq 0,$ the family $\{T^{j+Nk} \varphi\}_{k=0}^\infty$ is a frame for $\h$ for any $N \in \mn$ and any
$j\in \{0,1, \dots, N-1\}.$ \eni 
\et

\bp  As the first step of the proof we will show that
the assumptions imply that  $\{\lambda_k^N\}_{k=1}^\infty$ satisfies the Carleson condition for
all $N\in \mn.$ First, let $N=1$. By assumption the sequence $\{|\lambda_k|\}_{k=1}^\infty$
consists of distinct numbers, so Proposition  \ref{1406a} (ii) implies that
(\ref{1606c}) holds. Using now Proposition  \ref{1406a} (i) on the sequence $\{\lambda_k\}_{k=1}^\infty$, it follows that  $\{\lambda_k\}_{k=1}^\infty$ also
satisfies the Carleson condition, as claimed. Now, fix an arbitrary $N\in\mn$, $N>1$. Then
\bes \frac{1- |\lambda_{k+1}^N|}{1-|\lk^N|} = \frac{1-|\lambda_{k+1}|}{1-|\lk|}
\cdot
\frac{1+|\lambda_{k+1}|+\ldots +|\lambda_{k+1}|^{N-1}}{1+|\lambda_k|+\ldots +|\lambda_k|^{N-1}}.\ens
As in the first part of the proof,  take $c$ which satisfies (\ref{1606c}) and choose any
$\varepsilon\in (0, \frac{1-c}{c}).$
Using Proposition \ref{1406a} (iii),  we see  that
$$\lim_{k\to\infty}\frac{1+|\lambda_{k+1}|+\ldots +|\lambda_{k+1}|^{N-1}}{1+|\lambda_k|+\ldots +|\lambda_k|^{N-1}}=1
$$
and therefore there is $K_0\in\mn$ so that
$$\frac{1+|\lambda_{k+1}|+\ldots +|\lambda_{k+1}|^{N-1}}{1+|\lambda_k|+\ldots +|\lambda_k|^{N-1}}
<1+\varepsilon \ \mbox{ for } k>K_0.
$$
Then
for $k>K_0$ we have that
$$\frac{1-|\lambda_{k+1}|^N}{1-|\lambda_k|^N}
< c(1+\varepsilon)<1.
$$
By Proposition \ref{1406a} (i), the sequence $\{\lambda_k^N\}_{k=K_0+1}^\infty$ satisfies the  Carleson condition and hence  $\{\lambda_k^N\}_{k=1}^\infty$ satisfies the Carleson condition by Proposition \ref{1406a} (iv).

By Theorem \ref{19408bf} we know that for the  function $\varphi$ in \eqref{230412p},  $\{T^k\varphi\}_{k=0}^\infty$
is a frame for $\h.$ Now consider $N\in \mn, \, N\ge 2.$
Then the bounded operator $T^N$ is given by $T^N e_k=
\lambda_k^Ne_k,\, k\in \mn.$ Since $\{\lambda_k^N\}_{k=1}^\infty$  satisfies the
Carleson condition, Theorem 
\ref{19408bf}  shows
that $\{T^{Nk}\widetilde{\varphi}\}_{k=0}^\infty$ is
a frame for $\h$ whenever $\widetilde{\varphi}= \suk \widetilde{m_k} \sqrt{1- |\lambda_k^N|^{2}} e_k$ for some coefficients $\widetilde{m_k}$ that are bounded below and above. 
We aim at taking 
\bes  \widetilde{m_k}= m_k \sqrt{   \frac{1- | \lk|^2}{1- | \lk^N|^2}};\ens
note that this choice actually implies that
$\widetilde{\varphi}=\varphi.$ Now,

\bes \widetilde{m_k}= m_k \sqrt{   \frac{1- | \lk|^2}{1- | \lk^N|^2}} & = & m_k
\sqrt{\frac{1-|\lk|}{1-|\lk|} \cdot \frac{1+ | \lk|}{1+ | \lk| + \cdots + | \lk|^{2N-1}}} \\
& = & m_k \sqrt{\frac{1+ | \lk|}{1+ | \lk| + \cdots +| \lk|^{2N-1}}}.
\ens Using that $0 \le |\lk | \le 1$ and that $C_1 \le |m_k| \le C_2,$ we see that 
$C_1 (2N)^{-1/2} \le | \widetilde{m_k}| \le C_2.$  Thus $\widetilde{m_k}$ is indeed bounded below and above, i.e.,
we have proved that $\{T^{Nk}\varphi\}_{k=0}^\infty$ is a frame as stated in (i).

In order to prove (ii), fix any $N,K\in \mn,$  and note that 
\bes \{ T^{Nk}\varphi\}_{k=K}^\infty = \{ T^{NK} T^{Nk}\varphi\}_{k=0}^\infty.\ens  Since $1\ge |\lambda_k| \ge| \lambda_1| >0$ for all $k\in \mn,$ the bounded operator $T$ is
surjective, and hence 
$T^{NK}$ is  bounded and surjective; thus the result in (ii) follows immediately from (i) and
the well-known fact that the image of a frame
under a bounded surjective operator itself forms a frame.
The result in (iii) follows from the same argument
applied to the operator $T^j.$  \ep


\bex \rm{As a concrete example of a sequence $\lkn$ satisfying all the assumptions in Theorem \ref{19408b} , one can take $\lkn=\{1-\alpha^{-k}\}_{k=1}^\infty$, for any $\alpha >1.$ } \ep \enx

Note that the assumptions in Theorem \ref{19408b} are on the absolute
values $| \lk|.$
The following example shows that  $\{\lambda_k\}_{k=1}^\infty$ satisfying
the Carleson condition does not imply that $\{\lambda_k^N\}_{k=1}^\infty$ satisfies
the Carleson condition for $N>1.$

\begin{ex} Consider a sequence $\{\mu_k\}_{k=1}^\infty$ that satisfies the Carleson
	condition.
	Take a number $q\in(0,1)$ such that $\mu_k \neq q$ and $\mu_k\neq -q$ for every $k\in\mn$.
	By Proposition \ref{1406a} (iv), the sequence $\{\lambda_k\}_{k=1}^\infty:=\{q,-q, \mu_1,\mu_2,\mu_3,\ldots\}$ also satisfies the Carleson
	condition, but the sequence $\{\lambda_k^2\}_{k=1}^\infty$ does not.
\end{ex}

Theorem \ref{19408b} shows that for any $N\in \mn,$  the
considered class of
Carleson frames  can be
considered as the union of $N$ frames, namely,

\bee \label{10918b} \{T^k \varphi\}_{k=0}^\infty = \bigcup_{j=0}^{N-1} \{T^{Nk+j} \varphi\}_{k=0}^\infty.\ene

We will now show if we assume that $\{\lambda_k\}_{k=1}^\infty\subset \mathbb{D}$
consist of positive scalars, then ``one can weave the elements from these
frames from a sufficiently large index" (see Theorem \ref{10918f} for a precise statement). For this purpose we
need the subsequent lemma.

\bl  \label{10918c} Let $ \{T^k \varphi\}_{k=0}^\infty$ be
a Carleson frame as in Theorem \ref{carleson}, for which the sequence  $\{\lambda_k\}_{k=1}^\infty\subset \mathbb{D}$
consists of positive and strictly increasing scalars, and
fix  any $N\in \mn.$ Then, for any choice of
$j_k \in \{0, 1, \dots, N-1\}, \, k\in \mno,$
\bee \label{10918d} \sum_{k=0}^\infty || T^{Nk}\varphi - T^{Nk+ j_k}\varphi||^2 < \infty.\ene
\el

\bp By direct calculation, and using the upper bound
on the sequence $\{m_k\}_{k=0}^\infty$ in \eqref{232111a}

\bes \sum_{k=0}^\infty || T^{Nk}\varphi - T^{Nk+ j_k}\varphi||^2 & = &
\sum_{k=0}^\infty \sum_{n=1}^\infty  \left| m_n \lambda_n^{Nk} \sqrt{1-\lambda_n^2}- m_n \lambda_n^{Nk+j_k} \sqrt{1-\lambda_n^2} \right|^2 \\  & \le C_2^2 & \sum_{n=1}^\infty (1-\lambda_n^2) \sum_{k=0}^\infty \left(  \lambda_n^{2N}\right)^k \left(1- \lambda_n^{j_k}\right)^2. \ens
Using  that $0< 1-\lambda_n^{j_k} \le 1- \lambda_n^N \le 1-\lambda_n^{2N},$ this implies that
\bes \sum_{k=0}^\infty || T^{kN}\varphi - T^{kN+ j_k}\varphi||^2 & \le & C_2^2
\sum_{n=1}^\infty (1-\lambda_n^2) (1-\lambda_n^N)^2 \frac1{1-\lambda_n^{2N}} \\ & \le &
C_2^2 \sum_{n=1}^\infty (1-\lambda_n^2) \le \frac{C_2^2}{C_1^2} || \varphi||^2 < \infty,
\ens as desired.  \ep

\bt  \label{10918f}  Let $ \{T^k \varphi\}_{k=0}^\infty$ be
a Carleson frame as in Theorem \ref{carleson}, for which the sequence  $\{\lambda_k\}_{k=1}^\infty\subset \mathbb{D}$
consists of positive and strictly increasing scalars, and fix  any $N\in \mn.$ Then, for any choice of
$j_k \in \{0, 1, \dots, N-1\}, \, k\in \mno,$ there exists $J\in \mno$ such that the set
\bee \label{10918e} \{T^{Nk} \varphi\}_{k=0}^{J-1} \cup \{T^{Nk+j_k} \varphi\}_{k=J}^\infty  \ene
is a frame. \et

\bp Fix $N\in \mn;$ by Theorem \ref{19408b} we know that
$ \{T^{Nk} \varphi\}_{k=0}^\infty$ is a frame for $\h.$ Let $A$ denote
a lower frame bound, and  fix any
choice of
$j_k \in \{0, 1, \dots, J-1\};$ then, by Lemma \ref{10918c},
we can choose $J\in \mno$ such that
\bes \sum_{k=J}^\infty || T^{Nk}\varphi - T^{Nk+ j_k}\varphi||^2 < A.\ens
By standard perturbation results
in frame theory (see, e.g., page 565 in \cite{CB}),
this implies that in the frame $\{T^{Nk}\varphi\}_{k=0}^\infty$ we can replace
the elements $\{T^{Nk}\varphi\}_{k=J}^\infty$ by
$\{T^{Nk+j_k}\varphi\}_{k=J}^\infty$  while the frame
property is kept, as desired. \ep

Formulated in the terminology of weaving frames (see
\cite{CasBem}), Theorem  \ref{10918f} says that if we fix a sufficiently
high number $J$ of the first elements of the frame $ \{T^{Nk} \varphi\}_{k=0}^\infty,$ we can
weave the ``tail frames" $\{T^{Nk+j}\}_{k=J+1}^\infty, \, j=0, \dots, N-1,$  and
keep the frame property.

We conjecture that the very particular properties of the
frames in  Theorem \ref{19408b} imply that indeed
the frames  $\{T^{Nk+j}\}_{k=0}^\infty, \, j=0, \dots, N-1,$
are weaving for any $N=2,3,....;$ that is,
we conjecture that Theorem \ref{10918f} holds
with $J=0.$
Note that for $N=2,$
this would mean that the frames $\{T^{2k}\varphi\}_{k=0}^\infty$ and $\{TT^{2k}\varphi\}_{k=0}^\infty$ are weaving. Conditions
for this to hold are provided in \cite{CasBem};
however, they are formulated in terms of the frame bounds,
and the available
frame bound estimates (see \cite{CHP}) are not strong
enough to give the desired conclusion.

The ``subsampling property" in Theorem \ref{19408b} (i)
is a very special property for the considered frames, and
has not been identified for any other frame in the literature. Let
us illustrate this by showing that it can not hold for a number of important classes of frames.

\bex \rm{ For $a,b>0$ and $g\in \ltr,$ let
	$E_{mb}T_{na}g(x):= e^{2\pi i mbx}g(x-na),$ and assume
	that $\mts$ is a Gabor frame for $\ltr.$  Since
a necessary condition for such a family to be a frame
is that $ab \le 1,$ it is clear that ``subsamples"
$\{ E_{mb}T_{nNa}g\}_{m,n\in \mz}$ or $\{ E_{mNb}T_{na}g\}_{m,n\in \mz}$ never can be frames for
arbitrary large values of $N\in \mn.$ } \ep \enx

\bex \rm{For $a>1, b>0,$ and $\psi\in \ltr,$ let
	$D_{a^j}T_{kb}\psi(x):= a^{j/2}\psi(a^j x-kb),$ and
	 consider the wavelet system $\{ D_{a^j}T_{kb}\psi\}_{j,k\in \mz}.$ In contrast to
	 the Gabor case, there exist wavelet frames for
	 arbitrary choices of the parameters $a>1, b>0,$
	 see \cite{Da2}. That is, assuming that  
  $\{ D_{a^j}T_{kb}\psi\}_{j,k\in \mz}$ is a frame in general does not exclude that  $\{ D_{a^{Nj}}T_{kb}\psi\}_{j,k\in \mz}$
  	can be a frame, even for large values of $N\in \mn.$ 
  	However, for the frame constructions in the literature, there is a connection between the choice of parameters $a,b$ and the corresponding function $\psi,$ i.e.,
  	changing the parameters in a wavelet frame will in general
  	also force a change in the generating function (see
  	\cite{CG} for an illustration of this).
  	
  	For the case of a {\it band limited} frame
  	 $\{ D_{a^{j}}T_{kb}\psi\}_{j,k\in \mz}$ we can exclude the possibility that  $\{ D_{a^{Nj}}T_{kb}\psi\}_{j,k\in \mz}$ is a frame
  	 for large values of $N\in \mn.$ To be precise, assume
  	 that  $\{ D_{a^{j}}T_{kb}\psi\}_{j,k\in \mz}$ is a frame with bounds $A,B;$ then it is well-known that
  	 \bee \label{230712a} bA \le \sum_{j\in \mz} | \widehat{\psi}(a^j \ga)|^2 \le bB, \, a.e. \, \ga \in \mr.\ene 
  	 Assume now that 
  	 $\supp\, \widehat{\psi} \subset [-D, -C] \cup [C,D]$
  	 for some $C,D>0.$  Then, considering $\ga \in [C/4, C/2],$ it is clear that for sufficiently large values of
  	 $N\in \mn,$ we have 
  	 $ \sum_{j\in \mz} | \widehat{\psi}(a^{Nj} \ga)|^2=0,$
  	 i.e., $\{ D_{a^{Nj}}T_{kb}\psi\}_{j,k\in \mz}$
  	 can not be a frame.  }  \ep \enx

 To illustrate the issue further, we finally prove that despite the fact that
``arbitrarily sparse"  subfamilies of the Carleson frames
in Theorem \ref{19408b} still form frames, no
frame can have the property that arbitrary infinite subfamilies keep the frame property.

\bpr \label{81206b} No frame $\ftk$ for an infinite-dimensional Hilbert space $\h$ has the property that every infinite subsequence $\{f_k\}_{k\in I}$
is a frame for $\h.$\epr

\bp Let $\ftk$ denote an arbitrary frame for $\h,$ and let $\{e_j\}_{j=1}^\infty$
denote an orthonormal basis for $\h.$ We will inductively construct an infinite family
$\{g_\ell\}_{\ell=1}^\infty \subset \ftk$ which is not a frame. First, since $\ftk$ is a Bessel sequence, choose
$N_1\in \mn$ such that
\bes \sum_{k=N_1}^\infty | \la e_1, f_k\ra|^2 \le 1,\ens
and let $g_1:=f_{N_1}.$ Thus, the lower frame bound (if it exists) for
$\{f_k\}_{k=N_1}^\infty$ is at most $1.$
Now, choose $j_1>1$ such that $|\la e_{j_1}, f_{N_1}\ra|^2 \le 2^{-2},$ choose $N_2>N_1$ such that
\bes \sum_{k=N_2}^\infty | \la e_{j_1}, f_k\ra|^2 \le 2^{-2},\ens
and let $g_2:=f_{N_2}.$ Then
\bes
|\la e_{j_1}, f_{N_1}\ra|^2+ \sum_{k=N_2}^\infty | \la e_{j_1}, f_k\ra|^2 \le 2^{-1},\ens
i.e., the lower frame bound (if it exists) for $\{f_{N_1}\} \cup \{f_k\}_{k=N_2}^\infty$
is at most $2^{-1}.$

Now, choose $j_2>j_1$ such that
$|\la e_{j_2}, f_{N_1}\ra|^2 +|\la e_{j_2}, f_{N_2}\ra|^2
\le 2^{-3},$ choose $N_3>N_2$ such that
\bes \sum_{k=N_3}^\infty | \la e_{j_2}, f_k\ra|^2 \le 2^{-3},\ens
and let $g_3:=f_{N_3}.$ Then
\bes
|\la e_{j_2}, f_{N_1}\ra|^2 +|\la e_{j_2}, f_{N_2}\ra|^2
+ \sum_{k=N_3}^\infty | \la e_{j_2}, f_k\ra|^2 \le 2^{-2},\ens
i.e., the lower frame bound (if it exists) for $\{f_{N_1}\} \cup
\{f_{N_2}\} \cup \{f_k\}_{k=N_3}^\infty$
is at most $2^{-2}.$

Inductively, having constructed $g_1, \dots, g_\ell$ for some $\ell\in \mn$ (as well as the scalars $j_1, \dots, j_{\ell-1}$ and $N_1, \dots, N_\ell$),  
choose $j_\ell>j_{\ell-1}$ such that
$|\la e_{j_\ell}, f_{N_1}\ra|^2 +|\la e_{j_\ell}, f_{N_2}\ra|^2+
\cdots + |\la e_{j_\ell}, f_{N_\ell}\ra|^2
\le 2^{-\ell -1},$ choose $N_{\ell+1}>N_\ell$ such that
\bes \sum_{k=N_{\ell+1}}^\infty | \la e_{j_\ell}, f_k\ra|^2 \le 2^{-\ell -1},\ens
and let $g_{\ell+1}:=f_{N_{\ell+1}}.$ Then
\bes
|\la e_{j_\ell}, f_{N_1}\ra|^2 +|\la e_{j_\ell}, f_{N_2}\ra|^2+
\cdots + |\la e_{j_\ell}, f_{N_\ell}\ra|^2
+ \sum_{k=N_{\ell+1}}^\infty | \la e_{j_\ell}, f_k\ra|^2 \le 2^{-\ell},\ens
i.e., the lower frame bound (if it exists) for $\{f_{N_1}\} \cup
\{f_{N_2}\} \cup \cdots \cup  \{f_{N_{\ell}}\}\cup \{f_k\}_{k=N_{\ell+1}}^\infty$
is at most $2^{-\ell}.$  By construction, $\{g_j\}_{j=1}^\infty$ is a subsequence
of $\{f_{N_1}\} \cup
\{f_{N_2}\} \cdots \cup  \{f_{N_{\ell}}\}\cup \{f_k\}_{k=N_{\ell+1}}^\infty$
for any $\ell\in \mn,$ and thus no strictly positive lower bound exists. \ep

\vn{\bf Acknowledgment:} F. Philipp was funded by the Carl Zeiss Foundation within the project {\em DeepTurb--Deep Learning in and from Turbulence}. He was further supported by the Free State of Thuringia and the German Federal Ministry of Education and Research (BMBF) within the project {\em THInKI--Th\"uringer Hochschulinitiative für KI im Studium}. D. Stoeva acknowledges support from the Austrian Science Fund (FWF) through Project P 35846-N ``Challenges in Frame Multiplier Theory'',  grant DOI 10.55776/P35846.

 {\bf \noindent Ole Christensen, 
 	Department of Applied Mathematics and Computer Science,
 	Technical University of Denmark,
 	Building 303,
 	2800 Lyngby, Denmark\\
 	Email: ochr@dtu.dk
 	
 	\vn Marzieh Hasannasab, 	Department of Applied Mathematics and Computer Science,
 	Technical University of Denmark,
 	Building 303,
 	2800 Lyngby, Denmark\\
 	Email:  mhas@dtu.dk

 	\vn Friedrich Philipp, Technische Universit\"at Ilmenau,
 	Institute of Mathematics,
 	Weimarer Straße 25,
 	D-98693 Ilmenau, Germany  \\
 	Email:    Friedrich.Philipp@tu-ilmenau.de

\vn Diana Stoeva, Faculty of Mathematics, University of Vienna,
Oskar-Morgenstern-Platz 1, Vienna 1090, Austria \\
E-mail address: diana.stoeva@univie.ac.at

}
\end{document}